\documentclass[a4paper,leqno]{article}
\usepackage{amssymb,amsmath,amsthm}
\usepackage{graphicx}
\usepackage[all]{xy}
\usepackage[french]{babel}
\newcounter{nequa}
\newcommand{\equa}{\stepcounter{nequa}\thenequa )\ \ }
\def\e{\varepsilon}
\def\f{\varphi{}}
\def\O{\Omega}
\def\ens#1{\left\{ #1\right\}}

\def\p{{\mathbb P}^3}
\def\qq{\quad ,\quad}
\def\G{{\cal G}}

\def\ps{{\cal P}{\cal S}}
\def\e{\varepsilon}
\def\f{\varphi}
\def\pt{\p\tilde\times\p}
\def\ft{F\tilde\times F}
\def\ov{\overline}
\def\Im{\hbox{Im}}
\def\triv{{\mathbb I}}
\newtheorem*{prop}{Proposition}
\newtheorem*{lemme}{Lemme}
\newtheorem*{remarque}{Remarque}
\newcommand{\demo}{\noindent{\it Preuve: }}
\renewcommand{\qed}{\par\rightline{\raise 3pt\hbox{{\it q.e.d.\qquad\qquad}}}\medskip}
\def\h{\hbox{Hom}}
\begin{document}
\title{Le calcul de Schubert selon Schubert}
\author{Felice Ronga}

\maketitle
\begin{abstract}
    We try to understand and justify Schubert calculus the way 
    Schubert did it.
\end{abstract}

\section{Introduction}
Dans son livre \cite{Schubert:1879} ``Kalk\"ul der abz\"ahlende Geometrie'',
paru en 1879, le Dr. 
Hermann C. H. Schubert a pr\'esent\'e une m\'ethode pour r\'esoudre des 
probl\`emes de g\'eom\'etrie \'enum\'erative, qu'il a appliqu\'ee avec 
succ\`es \`a un grand nombre de cas. 

Le 15\`eme probl\`eme de Hilbert proposait de donner des bases 
rigoureuses pour justifier les r\'esultats de Schubert; cela a 
\'et\'e fait dans une large mesure, par la th\'eorie de 
l'intersection (voir \cite{Kleiman:1976},\cite{Kleiman-Laksov:1972}, 
\cite{Fulton:1984}), et la plupart des r\'esultats 
des calculs de Schubert ont pu \^etre confirm\'es.

Notre propos ici est de comprendre et justifier la m\'ethode m\^eme 
que  Schubert a utilis\'ee. On retracera aussi ses calculs dans 
quelques cas simples, dans le but d'illustrer sa d\'emarche.
Dans  \cite{Kazarian:2003} on trouve un traitement tr\`es g\'en\'eral, 
bas\'e sur les multisingularit\'es d'applications, de probl\`emes de 
g\'eom\'etrie \'enum\'erative, qui recouvrent entre autres les 
r\'esultats de Schubert concernant les contacts de sous-espaces 
lin\'eaires avec une hypersurface d'un espace projectif.

Voici sommairement en quoi 
consiste la m\'ethode de Schubert. On distingue tout d'abord des \'el\'ements de 
base de l'espace projectif complexe~: points, plans, droites (dans 
l'ordre de complexit\'e des ensembles de ces \'el\'ements). On 
repr\'esente par des 
symboles, disons $x$, $y$, des conditions (en allemand~: {\it 
Bedingungen})\/ que doivent 
satisfaire des objets g\'eom\'etriques; le produit $x\cdot y$ de deux 
conditions repr\'esente la condition que $x$ et $y$ doivent \^etre 
satisfaites, la somme $x+y$ repr\'esente la condition que $x$ ou $y$ 
est satisfaite. Les conditions que l'on peut mettre sur
les \'el\'ements de base, \`a l'aide d'autre \'el\'ements de base 
(par exemple : les droites qui passent par un point) satisfont un 
certain nombre de formules, que l'on peut d\'eterminer assez 
simplement par des raisonnements g\'eom\'etriques, et des r\`egles de 
calcul \'el\'ementaire. 

Pour r\'esoudre un probl\`eme de 
g\'eom\'etrie \'enum\'erative, on s'efforce de l'exprimer en termes 
de conditions sur les \'el\'ements de base, en ayant \'eventuellement 
recours \`a des situations d\'eg\'en\'er\'ees, g\'eom\'etriquement 
plus simples \`a traiter, quitte \`a devoir tenir compte de 
multiplicit\'es des solutions. Schubert justifie  cette fa\c con de 
proc\'eder par le principe de la conservation du nombre ({\it Prinzip des Erhaltung der 
Anzahl}, \cite[\S\ 4, page 12]{Schubert:1879})\/, qui dit en gros 
que dans les cas 
d\'eg\'en\'er\'es, le nombre de solutions d'un probl\`eme 
\'enum\'eratif est conserv\'e, pourvu qu'il soit fini. Comme toute 
justification de ce principe, Schubert se base sur l'analogue 
alg\'ebrique~: le nombre de solutions d'une \'equation (polyn\^omiale) 
ne change pas lorsqu'on en fait varier les coefficients, sauf si elle 
devient une identit\'e, auquel cas on a une infinit\'e de solutions.

La force de la d\'emarche de Schubert r\'eside dans le fait que la 
notation symbolique, et les op\'erations sur les symboles, contiennent en 
germe la notion d'anneau de cohomologie (ou si l'on pr\'ef\`ere, 
l'anneau de Chow). Une condition $x$ repr\'esente en fait 
une famille de conditions, que l'on peut interpr\'eter comme une 
classe cohomologie d'un espace de configurations. Il est sous-entendu que lorsqu'on fait 
un produit, disons $x\cdot y$, les ensembles des objets auxquels 
s'adressent ces conditions doivent \^etre en position g\'en\'erale 
(ou tout au moins leur intersection doit avoir la bonne dimension, quitte \`a invoquer une 
multiplicit\'e). Un exemple~: $p$ 
d\'esigne la condition, adress\'ee aux points de ${\mathbb P}^3$, 
d'\^etre contenus dans un plan; $p_{g}$ d\'esigne la condition que 
les points soient sur la droite $g$; on a alors~: $p\cdot p=p_{g}$. En 
d'autres termes, le plan qui exprime la condition $p$ est g\'en\'erique; 
sans cela on aurait $p\cdot p=p$. En somme, l'ambigu\"\i{}t\'e 
\'eclair\'ee du 
calcul avec les symboles 
est ce qui fait son efficacit\'e. Il faut relever qu'avec quelques 
bons principes et des raisonnements g\'eom\'etriques simples, Schubert 
est parvenu \`a obtenir une quantit\'e de r\'esultats remarquables, 
pour la justification desquels, selon les crit\`eres de rigueur 
actuels, des d\'eveloppements importants ont \'et\'e n\'ecessaires.

En plus des formules sur les \'el\'ements de base, Schubert a  \'etabli
ce qu'il a appell\'e {\it Coinzidenzformeln}\/, principalement la 
formule $1)$ de la page 44 de \cite{Schubert:1879}, 
qui est un prototype de 
formule d'intersection r\'esiduelle, telle qu'on la trouve par exemple dans 
\cite[th\'eor\`eme 9.2]{Fulton:1984}.
Elle lui a permis entre autres d'\'etablir des formules de 
co\"\i{}ncidences multiples ({\it mehrfache Coinzidenzen}\/), avec une 
efficacit\'e et une rigueur qui n'ont rien \`a envier \`a leurs 
versions plus r\'ecentes que l'on trouve par exemple dans \cite{Kazarian:2003}.

En terme de cohomologie, si $X$ est un espace de configuration 
d'objets g\'eom\'etriques (par exemple les points d'une surface, l'espace des 
coniques), une condition $x$ peut \^etre repr\'esent\'ee comme la classe duale 
\`a une classe d'homologie de cycles $[\Omega_{x}]$ de $X$. Alors la condition $x\cdot y$ 
repr\'esente la clase duale \`a $\Omega_{x}\cap \Omega_{y}$, pourvu 
que ces deux cycles soient en position g\'en\'erale. Les formules que 
Schubert d\'emontre pour les \'el\'ements de base correspondent aux 
calculs des anneaux de cohomologie de l'espace projectif complexe 
${\mathbb  P}^3$, de son dual $\check{\mathbb  P}^3$, 
de la grassmannienne ${\cal G}$ des droites de  
${\mathbb  P}^3$, et finalement de l'espace ${\cal PS}$ (Punkt und 
Strahl), constitu\'e des paires form\'ees d'une droite de l'espace 
projectif et d'un point 
sur la droite.

Se donner une droite de $\p$ revient \`a se donner un 
sous-espace vectoriel de dimension $2$ de ${\mathbb C}^4$, et donc on peut aussi 
regarder la grassmannienne $\G$ 
comme l'ensemble des sous-espaces vectoriels de dimension $2$ de 
${\mathbb C}^4$; on voit alors qu'on a un fibr\'e vectoriel de rang 
$2$ sur $\G$, appel\'e fibr\'e tautologique, qui est constitu\'e des 
paires $(\alpha,v)$, o\`u $\alpha\in\G$ et $v\in\alpha$. En fait, 
$\ps$ n'est autre que le fibr\'e projectif associ\'e 
\`a $\eta$. 

Notons que Schubert n'avait pas introduit de 
notation pour d\'esigner les espaces ${\mathbb  P}^3$, $\check{\mathbb  P}^3$, 
${\cal G}$ et ${\cal PS}$,
puisqu'ils consituent en 
quelque sorte l'univers ambiant.

\section{Formules pour les espaces de configurations de 
base}
On introduit des symboles qui d\'esignent des objets 
g\'eom\'etriques dans les divers espaces de base ${\mathbb  P}^3$, $\check{\mathbb  P}^3$, 
${\cal G}$ ou ${\cal PS}$. Ces m\^emes symboles d\'esigneront aussi 
des conditions impos\'ees aux objets de base. Les ensembles ainsi 
d\'efinis engendrent l'homologie des espaces de base; dans le cas de 
$\p$ et $\G$ ils forment m\^eme une d\'ecomposition cellulaire 
minimale. 
En exprimant leurs intersections en termes des 
\'el\'ements de base, on aura d\'etermin\'e l'anneau de cohomologie de 
ces espaces.

Nous utiliserons \'evidemment la m\^eme notation que Schubert, bas\'ee 
sur les noms allemands des divers objets. Aussi est-il utile de 
rappeler quelques mots de la langue allemande~:
$$
\begin{tabular}{ll}
    Punkt :& point\\
    Gerade :& droite\\
    Ebene :& plan\\
    Strahl:& litt\'eralement : rayon; ici d\'esigne le plus souvent l'ensemble des droites dans un plan\\
	   & fix\'e, passant par un point fix\'e, c'est-\`a-dire un pinceau de droites.\\
	   & Parfois, ce mot est 
	   synonyme
	   de droite, comme dans ``Punkt und Strahl''
\end{tabular}
$$

Notons enfin que, \`a d\'efaut d'un mot plus pr\'ecis, nous utiliserons le mot {\it condition}, traduction du 
mot allemand {\it Bedingung}\/, pour signifier une exigence \`a 
laquelle des objets g\'eom\'etriques sont astreints.

Nous travaillerons avec l'anneau de cohomologie, mais on peut tout 
aussi bien utiliser l'anneau de Chow.

Lorsque des formules sont num\'erot\'ees, leur num\'ero est le m\^eme 
que dans \cite{Schubert:1879}.

\subsection{L'espace projectif complexe ${\mathbb  P}^3$}
Voici les conditions de base (ou simples) que l'on peut mettre sur 
les points de l'espace~:
$$
\begin{tabular}{cl}
    Notation &\hfil Condition \hfil\\
    $p$ & le point doit \^etre dans un plan donn\'e\\
    $p_g$ & le point doit \^etre sur une droite donn\'ee\\
    $P$ & le point lui-m\^eme est donn\'e
    \end{tabular}
$$
On v\'erifie facilement les relations~:
$$
\equa
p^2=p_g
\qq
\equa p^3=p\cdot p_g
\qq
\equa p\cdot p_{g}=P
\qq
\equa p^3=P
\quad .
$$
A titre d'exemple, l'interpr\'etation g\'eom\'etrique p\'edante de la premi\`ere formule 
est la suivante~: soient $e_1 ,e_2\subset {\mathbb  P}^3$ deux plans et
$$
\Omega_{e_i}=\{P\in{\mathbb  P}^3\mid P\in e_i\}
\quad ,\quad
i=1,2
\quad .
$$
Alors $p^2$ d\'esigne les points de l'intersection $\Omega_{e_1}\cap\Omega_{e_2}$ 
lorsque $e_1$ et $e_2$ sont en position g\'en\'erale, c'est-\`a-dire 
lorsque $e_1$ et $e_2$ se coupent en une droite; l'ensemble des points 
astreints \`a se trouver sur cette droite est not\'e $p_g$.

Voici maintenant l'interpr\'etation en cohomologie de ces formules. 
D\'esignons par $t\in H^2({\mathbb  P}^3,{\mathbb Z})$ la classe 
duale au cycle consitu\'e par les points d'un plan de ${\mathbb  P}^3$; 
alors $t^2$ est la classe duale \`a une droite et $t^3$ est la classe 
duale \`a un point de ${\mathbb  P}^3$.

Notons que si on choisit un drapeau $p\in g\subset e$, en notant 
par
$\Omega_p$, $\O_g$ et $\O_e$ les ensembles de points correspondants, 
on a que $\O_p\subset\O_g\subset\O_e\subset\p$ est une 
d\'ecomposition cellulaire de $\p$.

On peut traiter de mani\`ere analogue le cas de $\check\p$, l'espace 
des plans de $\p$~:
$$
\begin{tabular}{cl}
    Notation&\hfil Condition\hfil \\
    $e$ & le plan doit passer par un point donn\'e\\
    $e_g$ & le plan doit contenir une droite donn\'ee\\
    $E$ & le plan lui-m\^eme est donn\'e
    \end{tabular}
$$
On a les formules~:
$$
\equa e^2=e_g
\qq
\equa e^3=e\cdot e_{g}
\qq
\equa e\cdot e_g =E
\qq
\equa e^3=E
\quad .
$$
\subsection{La grassmannienne ${\cal G}$ des droites de $\p$}
Voici les conditions de base~:
$$
\begin{tabular}{clc}
    Notation &\hfil Condition\hfil &Dimension \\
    $g$ & la droite doit  couper une droite donn\'ee&3\\
    $g_e$ & la droite doit se trouver dans un plan donn\'e&2\\
    $g_p$&la droite doit passer par un point donn\'e&2\\
    $g_s$ & la droite doit appartenir \`a un pinceau donn\'e&1\\
    $G$& la droite elle-m\^eme est donn\'ee&0
    \end{tabular}
$$
Si on choisit un drapeau $P\in g\subset e\subset\p$, notons 
 $\O_{g}$, $\O_{e}$, $\O_{p}$, $\O_{s}$, $\O_{G}=\ens{G}$ 
les droites satisfaisant respectivement les conditions $g$, $g_{e}$, 
$g_{p}$, $g_{s}$, $G$. On a un diagramme d'inclusions~:
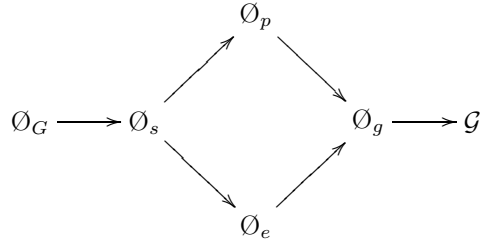
$$
\xymatrix{
&&\O_{p}\ar[dr]&&\\
\O_G\ar[r]&\O_{s}\ar[ur]\ar[dr]&&\O_{g}\ar[r]&{\cal G}\\
&&\O_{e}\ar[ur]&&
}
$$
et les $\O_{.}$ sont les cellules d'une d\'ecomposition cellulaire 
de ${\cal G}$ (voir \cite[\S\ 6]{Milnor:1974}). On appelle ces cellules {\it cycles de Schubert}\/.

Pour exprimer $g^2$ en termes des conditions de base, on suppose que 
les deux droites donn\'ees $g$ et $g'$ se coupent en un point $P$; en 
prenant pour $e$ le plan de $g$ et $g'$ on a~:
$$
\O_{g}\cap\O_{g'}=\O_{p}\cup \O_{e}
$$
et de l\`a Schubert conclut, en invoquant le principe de la 
conservation du nombre, la formule~:
$$
\equa g^2=g_{p}+g_{e}\quad .
$$
On va justifier cette formule de deux mani\`eres~: d'abord, en 
exprimant les calculs en cohomologie. Ensuite, en montrant que 
$\O_{g}$ et $\O_{g'}$ se coupent transversalement en dehors du lieu 
$\O_{s}$ des droites du plan $e$ passant par $P$, qui est de dimension 
$1$; ceci justifie en fait la d\'emarche de Schubert~: malgr\'e que 
la situation o\`u $g$ et $g'$ sont coplanaires soit 
d\'eg\'en\'er\'ee, il n'y a pas de multiplicit\'e \`a prendre en 
compte dans l'intersection.

Notons que le groupe lin\'eaire ${\cal G}\ell(4,{\mathbb C})$ agit 
transitivement sur ${\cal G}$. Il suit de l\`a que si l'on prend des 
droites, des points ou des plans g\'en\'eriques, les cycles de Schubert 
correspondants sont transverses. On v\'erifie alors facilement 
les formules suivantes~:
$$
\displaylines{
    \equa g\cdot g_p=g_s\qq\equa g\cdot g_e=g_s\cr
    \equa g\cdot g_s=G\qq\equa g_p\cdot g_e=0
    }
$$
En multipliant la formule $9)$ par $g$ et en utilisant $10)$ et $11)$ on obtient~:
$$
\equa g^3= g\cdot g_p+g\cdot g_e
\qq
\equa g^3=2\cdot g_s
$$
En multipliant encore par $g$~:
$$
\equa g^4=2\cdot g\cdot g_s=2\cdot g^2\cdot g_e=2\cdot g^2\cdot g_p =
2\cdot g_p^2=2\cdot g_e^2=2\cdot G
\quad .
$$
Notons que la formule $g^4=2\cdot G$ nous dit que le nombre de 
droites qui s'appuyent sur $4$ droites donn\'ees est \'egal \`a $2$. 
C'est un premier exemple, souvent cit\'e, d'application du calcul de Schubert \`a la 
g\'eom\'etrie \'enum\'erative.

\subsubsection{Cohomologie de ${\cal G}$}

Regardons la grassmannienne $\G$ comme l'espace des sous-espaces vectoriels de 
dimension $2$ de ${\mathbb C}^4$.
Soit $\eta=(E\stackrel{\pi}{\to} \G)$ le fibr\'e tautologique,
fibr\'e vectoriel de rang $2$~:
$$
E=\ens{(\alpha,v)\in\G\times{\mathbb C}^4\mid v\in \alpha}
\qq
\pi(\alpha,v)=\alpha
\quad .
$$
Soient $c_{i}(\eta)\in H^{2i}(\G,{\mathbb Z})$, $i=1,2$, les classes de 
Chern de $\eta$, et $s_{i}(\eta)\in H^{2i}(\G,{\mathbb Z})$, 
$i=1,\dots ,4$, les classes de Segre (voir par exemple 
\cite{Fulton:1984}). 
Elles sont li\'ees par la relation~:
$$
(1+c_{1}(\eta)+c_{2}(\eta))\cdot(1+s_{1}(\eta)+s_{2}(\eta)+s_{3}(\eta)
+s_{4}(\eta))=1
\quad .
$$
D\'esignons par $\triv^n$ le fibr\'e trivial de rang $n$, de base 
un espace non sp\'ecifi\'e qu'on peut d\'eduire du contexte. Puisque 
$\eta\subset\triv^4$, on peut poser $\eta'=\triv^4/\eta$; alors 
$c(\eta')=s(\eta)$.

Soient maintenant $x_{1}$ et $x_{2}$ des variables formelles et soient
$y_{1},y_{2}\in{\mathbb Z}[x_{1},x_{2}]$ d\'efinis par la relation~:
$$
(1+x_{1}+x_{2})\cdot(1+y_{1}+y_{2}+y_{3}+y_{4})=1
$$
ce qui revient \`a poser~:
$$
y_{1}=-x_{1}
\qq
y_{2}=x_{1}^2-x_{2}
\qq
y_{3}=2x_{1}x_{2}-x_{1}^3
\qq
y_{4}=x_{1}^4-x_{2}^2+3x_{1}^2x_{2}
$$
comme on v\'erifie facilement. On peut montrer (voir \cite[proposition 
page 69]{Stong:1968}) que l'homomorphisme 
d'anneaux~:
$$
{\mathbb Z}[x_{1},x_{2}]\to H^*(\G,{\mathbb Z})
\qq
x_{i}\mapsto c_{i}(\eta)
$$
induit un isomorphisme d'anneaux~:
$$
{\mathbb Z}[x_{1},x_{2}]/I(y_{3},y_{4})
\stackrel{\simeq}{\longrightarrow} H^*(\G,{\mathbb Z})
$$
 o\`u $I(y_{3},y_{4})$ d\'esigne l'id\'eal engendr\'e par 
 $y_{3}$ et $y_{4}$. On voit que $H^*(\G,{\mathbb Z})$ est engendr\'e en 
 tant que groupe par~:
 $$
 c_{1}
 \qq 
 c_{1}^2
 \qq
 c_{2}
 \qq
 c_{1}c_{2}
 \qq
 c_{2}^2
 $$
 et la structure d'anneau est determin\'ee par les relations 
 $2c_{1}c_{2}-c_{1}^3=0$, $c_{1}^4-c_{2}^2+3c_{1}^2c_{2}=0$, d'o\`u on 
 d\'eduit encore que $2c_{1}^2c_{2}-c_{1}^4=0$ et 
 $c_{1}^2c_{2}=c_{2}^2$. 
 \begin{remarque}
     Dans \cite[proposition page 69]{Stong:1968}, on affirme que 
     $H^*(\G)\simeq{\mathbb Z}[c_{1},c_{2}]/I(\ens{s_{j},j>2})$ o\`u les 
     $s_{j}$ sont d\'efinies pour tout entier positif
     par la relation~:
     $$
     (1+c_{1}+c_{2})(1+s_{1}+s_{2}+\cdots +s_{j}+\cdots)=1
     $$
     \`a valoir dans l'anneau gradu\'e ${\mathbb Z}[c_{1},c_{2}]$.
     Mais on voit facilement que $s_{j}\in I(s_{1},\dots ,s_{j-1})$, et 
     donc
     $$
     I(\ens{s_{j},j>2})=I(s_{3},s_{4})\quad .
     $$
 \end{remarque}
 On va exprimer les classes duales par Poincar\'e des diverses 
 cellules en termes des classes de Chern et de Segre de $\eta$. Voici 
 d\'ej\`a le r\'esultat~:
 $$
 \begin{tabular}{|c|c|c|c|c|c|c|}\hline
     Notation symbolique & -- & $g$ & $p$ & $e$ & $s$ & $G$ \\ \hline
     Cycle & $\G$ & $\O_{g}$ & $\O_{p}$ & $\O_{e}$ & $\O_{s}$ & $\O_{G}$ \\ \hline
      Classe duale & 1 &$s_{1}$ & $s_{2}$ & $c_{2}$ & $s_{1}c_{2}$ & 
      \vphantom{$\int_{0}^1$}$c_{2}^2=s_{2}^2$ \\ 
     \hline
 \end{tabular}
 $$

Pour cela, appelons $v_{i}$, $i=1,\dots ,4$ une base de ${\mathbb C}^4$; on 
d\'esignera par $\langle v_{i_{1}},\dots ,v_{i_{k}}\rangle$ l'espace 
vectoriel engendr\'e par $v_{i_{1}},\dots ,v_{i_{k}}$. On va exprimer 
les diverses conditions
avec les \'el\'ements du drapeau
$$
P=\langle v_{1}\rangle\subset g=\langle v_{1},v_{2}\rangle\subset
e=\langle v_{1},v_{2},v_{3}\rangle\subset{\mathbb C}^4
$$
\noindent\fbox{$\O_{g}$}

Consid\'erons le morphisme de fibr\'es $\varphi_{g}:\eta\to\triv^4/\langle 
v_{1},v_{2}\rangle$ induit par l'inclusion naturelle de $\eta$ dans 
$\triv^4$. En se rappelant que $g$ est la droite qui est le 
projectivis\'e de l'espace vectoriel $\langle v_{1},v_{2}\rangle$,
on voit que $\O_{g}=\Sigma(\varphi_{g})$, et donc sa classe duale est 
$s_{1}(\eta)$.

\vspace{\medskipamount}
\noindent
\fbox{$\O_{e}$} 

On consid\`ere le morphisme naturel $\varphi_{e}:\eta\to\triv^4/\langle 
v_{1},v_{2},v_{3}\rangle$, qui \'equivaut \`a une section $\sigma$ de
$\eta^*\otimes\triv^4/\langle v_{1},v_{2},v_{3}\rangle$, et $\O_{e}$ 
est \'egal aux z\'eros de cette section, ce qui fait que sa classe 
duale vaut $c_{2}(\eta^*)=c_{2}(\eta)$.

\vspace{\medskipamount}
\noindent
\fbox{$\O_{p}$} 

Ici on prend le morphisme naturel $\varphi_{p}:\langle 
v_{1}\rangle\to\triv^4/\eta$, qu'on peut voir comme section de 
$\triv^4/\eta$; ses z\'eros constituent $\O_{p}$, donc sa classe duale 
vaut $c_{2}(\triv^4/\eta)=s_{2}(\eta)$.

\vspace{\medskipamount}
\noindent
\fbox{$\O_{s}$ et $\O_{G}$}

Soient $e'$ le plan projectif correspondant \`a $\langle 
v_{1},v_{2},v_{4}\rangle$ et $g'$ la droite projective correspondante 
\`a $\langle v_{1},v_{4}\rangle$.
On remarque que $\O_{s}=\O_{g'}\cap\O_{e}$ et 
$\O_{G}=\O_{e}\cap\O_{e'}$, les intersections \'etant transverses. On 
en d\'eduit que les classes duales sont respectivement 
$s_{1}c_{2}$ et $c_{2}^2$.

\vspace{\medskipamount}
Par exemple, on peut retrouver la formule $9)$ en remarquant que 
 $s_{1}^2=c_{1}^2=(c_{1}^2-c_{2})+c_{2}=s_{2}+c_{2}$.
 
Aussi 
$s_{1}^4=s_{1}(-c_{1}^3)=s_{1}(-2c_{1}c_{2})=2c_{1}^2c_{2}=2c_{2}^2$ 
nous red\'emontre que $g^4=2G$.

On peut retrouver de fa\c  con analogue les autres formules.
\subsubsection{Justification par le principe de la conservation du nombre}

Pour introduire des coordonn\'ees locales sur $\G$, on choisit un 
sous-espace vectoriel  $\alpha_{0}\subset{\mathbb C}^4$ de dimension $2$
et un espace vectoriel suppl\'ementaire $\alpha'$. On d\'esigne par 
$\h(\alpha_{0},\alpha')$ l'espace des applications lin\'eaires de 
$\alpha_{0}$ dans $\alpha'$; on d\'efinit 
$\f:\h(\alpha_{0},\alpha')\to\G$ en associant \`a $A\in 
\h(\alpha_{0},\alpha')$ son graphe, ce qui d\'efinit une bijection 
sur l'ouvert
$$
U_{\alpha_{0},\alpha'}=\ens{\beta\in\G\mid \beta\cap\alpha'=\ens{0}}
\quad .
$$
On v\'erifie que l'on d\'efinit ainsi un atlas sur $\G$; on notera 
par $\ell_{A}$ la droite projective correspondante \`a $A\in 
\h(\alpha_{0},\alpha')$.
\begin{lemme}
    Soient $A,B\in \h(\alpha_{0},\alpha')$ et supposons que 
    $\ell_{A}\in\O_{\ell_{B}}$. Alors il existe une droite 
    vectorielle $\ell_{0}\subset\alpha_{0}$ telle que 
    $A|\ell_{0}=B|\ell_{0}$. 
    
    $\ell_{A}$ est un point r\'egulier de $\O_{\ell_{B}}$ si et 
    seulement si $A\not= B$, et si c'est le cas, on a~:
    $$
    T(\O_{\ell_{B}})_{\ell_{A}}
    =\ens{\ov A\in \h(\alpha_{0},\alpha')\bigm|
    \;\ov A|\ell_{0}:\ell_{0}\to\alpha'/\Im(A-B)\hbox{ est nulle }}
    $$
\end{lemme}
\demo
Au lieu de d\'ecrire $\O_{\ell_{B}}$ au voisinage de $\ell_{A}$, il 
est plus commode de se placer dans l'espace 
$\h(\alpha_{0},\alpha')\times \h(\ell_{0},\ell')$, o\`u $\ell'$ est 
une droite vectorielle suppl\'ementaire \`a $\ell_{0}$ dans 
$\alpha_{0}$. D\'esignons par 
$i_{\ell_{0}}:\ell_{0}\subset\alpha_{0}$ l'inclusion, par
$p:\h(\alpha_{0},\alpha')\times \h(\ell_{0},\ell'\to \h(\alpha_{0},\alpha'))$ 
la projection; 
l'\'equation 
$$
(A'-B)\circ( i_{\ell_{0}})+\lambda=0
\qq
A'\in \h(\alpha_{0},\alpha')
\, ,\,
\lambda\in \h(\lambda_{0},\lambda')
$$
d\'efinit un sous-ensemble $\tilde\O$ qui est en bijection par $p$ 
avec $\Omega_{B}\cap U_{\alpha_{0},\alpha'}$, sauf au-dessus de 
$A'=B$. Si on d\'erive cette \'equation en $A'=A$ on trouve~:
$$
\ov A\circ i_{\ell_{0}}+(A-B)\circ\ov\lambda=0
\quad ;
$$
si $A\not=B$, $\hbox{Ker}(A-B)=\ell_{0}$ et alors
$$
\exists\; \ov\lambda\hbox{ tel que }\ov A\circ i_{\ell_{0}}+(A-B)=0
\ \Longleftrightarrow\
\ov A\circ i_{\ell_{0}}:\ell_{0}\to\alpha'/\Im(A-B)\hbox{ est nulle }
$$
\qed
\begin{prop}
    Soient $\ell_{B_{1}}$ et $\ell_{B_{2}}$ deux droites distinctes, 
    se coupant en un point $P_{1,2}$. Alors $\O_{\ell_{B_{1}}}$ et
    $\O_{\ell_{B_{2}}}$ se coupent transversalement en dehors de 
    l'ensemble des 
    droites par $P_{1,2}$ se trouvant dans le plan par $\ell_{B_{1}}$ 
    et $\ell_{B_{2}}$. 
\end{prop}
\demo
Soit $\ell_{A}\in\O_{\ell_{B_{1}}}\cap\O_{\ell_{B_{2}}}$. Supposons 
d'abord que $\ell_{A}$ passe par $P_{1,2}$, et donc ne soit pas dans 
le plan de $\ell_{B_{1}}$ et $\ell_{B_{2}}$. Soit 
$\ell_{1,2}\subset\alpha_{0}$ la droite vectorielle correspondante 
\`a $P_{1,2}$, c'est-\`a-dire telle que 
$B_{1}|\ell_{1,2}=B_{2}|\ell_{1,2}=A|\ell_{1,2}$. Il suit du lemme que
$$
T(\O_{\ell_{B_{1}}})_A\cap T(\O_{\ell_{B_{2}}})_A
=\ens{\ov A\;\bigm|\; \ov A|\ell_{1,2}:\ell_{1,2}\to 
\alpha'/\Im(A-B_{i})\hbox{ est nulle },\; i=1,2}
$$
Puisque $A-B_{1}$ et $A-B_{2}$ ont m\^eme noyau $\ell_{1,2}$, si elles 
avaient m\^eme image on aurait~:
$$
A-B_{1}=\lambda A-B_{2}
$$
o\`u $\lambda$ est un scalaire, et 
$\lambda\not=1$, sans quoi 
$B_{1}=B_{2}$. On en d\'eduirait que
$$
A=\frac{1}{1-\lambda}B_{1}-\frac{\lambda}{1-\lambda}B_{2}
$$
et donc $\ell_{A}$ serait dans le plan par $\ell_{B_{1}}$ et 
$\ell_{B_{2}}$, contradiction. On a donc bien que les deux conditions 
que $\ov A|\ell_{0}\to\alpha'/\Im(A-(B_{i})$ soient nulles pour $i=1,2$  
sont ind\'ependantes, et donc on a transversalit\'e.

Si $\ell_{A}$ est dans le plan de $\ell_{B_{1}}$ et $\ell_{B_{2}}$ 
mais ne passe pas par $P_{1,2}$, soient 
$P_{1}=\ell_{A}\cap\ell_{B_{1}}$ et $P_{2}=\ell_{A}\cap\ell_{B_{2}}$
et soient $\ell_{1}$, $\ell_{2}\subset\alpha_{0}$ les droites 
d\'etermin\'ees par~:
$$
(A-B_{1})|\ell_{1}=0\qq
(A-B_{2})|\ell_{2}=0
\quad .
$$
Alors~:
$$
T(\O_{\ell_{B_{1}}})_A\cap T(\O_{\ell_{B_{2}}})_A
=\ens{\ov A\;\bigm|\; \ov A|\ell_{i}:\ell_{i}\to 
\alpha'/\Im(A-B_{i})\hbox{ est nulle },\; i=1,2}
$$
et comme $\ell_{1}\not=\ell_{2}$, ces deux conditions sont 
ind\'ependantes et on a transversalit\'e.
\qed
Soient $P,Q,R,S\in\p$ $4$ points non coplanaires, dont $3$ d'entre eux 
ne sont jamais align\'es. Alors, si on prend les quatre droites 
$\ell_{P,Q}$ par $P$ et $Q$, $\ell_{Q,R}$, $\ell_{R,S}$ et 
$\ell_{S,P}$, on voit que les cycles de Schubert correspondant se coupent 
transversalement en les deux droites $\ell_{P,R}$ et en $\ell_{Q,S}$. 
En effet, les intersections deux \`a deux sont transverses par la 
proposition 1, et la transversalit\'e des intersections restantes est 
\'el\'ementaire \`a v\'erifier (par exemple : intersection de 
l'ensemble des droites dans le plan par $P,Q,S$ et l'ensemble des 
droites dans $Q,R,S$).

Le probl\`eme \'enum\'eratif de trouver le nombre de droites 
s'appuyant sur quatre droites donn\'ees est souvent cit\'e en exemple pour illustrer 
les m\'ethodes de Schubert (voir \cite{Kleiman-Laksov:1972}).
Ce qui pr\'ec\`ede permet de justifier le recours
au cas un peu 
d\'eg\'en\'er\'e, o\`u les $4$ droites donn\'ees se coupent par 
paires. Par contre, comme Alexandre Gabard me l'a fait remarquer, si 
les quatre droites donn\'ees sont les droites du m\^eme syst\`eme 
d'une quadrique lisse, toute droite de l'autre syst\`eme rencontre 
ces quatre droites, et il y a donc une infinit\'e de solutions~: les 
4 droites donn\'ees ne sont pas en position g\'en\'erale du point de 
vue de ce probl\`eme.

Pour voir un cas d\'eg\'en\'er\'e ayant une seule solution avec 
multiplicit\'e $2$, on peut prendre pour $\ell_{1}$, $\ell_{2}$ et 
$\ell_{3}$ trois droites du m\^eme syst\`eme sur une quadrique lisse, 
pour $\ell_{4}$ une droite tangente en un point $P$ \`a la quadrique (mais pas 
contenue dedans). La seule solution est alors la droite par $P$ qui 
appartient \`a l'autre syst\`eme de droites.

\subsection{L'espace $\ps$ des points sur les droites de $\p$}
Rappelons que l'espace $\ps$ est consitu\'e des couples form\'es par 
un point sur une droite de l'espace. Pour exprimer des conditions,
on va utiliser des symboles de la forme $xy$, o\`u 
$x$ est un symbole qui exprime une condition sur les points, $y$ 
sur les droites. Ainsi, le symbole $pg$ d\'esigne les couples 
form\'es d'un 
point et une droite, le point \'etant astreint \`a se trouver dans un 
plan, la droite devant s'appuyer sur une droite donn\'ee; si on note 
$\O_{pg}$ l'ensemble de ces couples, et encore $g$ une droite et $e$ 
un plan
provisoirement fix\'es ~:
$$
\O_{pg}=\ens{(\ell,Q)\in\ps\mid \ell\cap g\not=\emptyset\qq Q\in e}
\quad .
$$
Remarquons que
$$
\O_{pg}=\ens{(\ell,Q)\in\ps\mid Q\in 
g}\cup\ens{(\ell,Q)\in\ps\mid\ell\subset e}
$$
d'o\`u on d\'eduit la formule (\cite[page 25]{Schubert:1879})~:
$$
\hbox{I})\ \ pg=p_{g}+g_{e}=p^2+g_{e}
\quad .
$$
Cette formule est fondamentale, dans le sens que toute autre formule 
est cons\'equence de celle-ci et des formules d\'ej\`a obtenues dans 
$\G$ et $\p$; la raison en est expliqu\'ee au \S\ suivant. Etablissons 
tout de m\^eme les formules suivantes; en multipliant $\hbox{I})$ par $p$, puis 
par $g$ on obtient~:
$$
\displaylines{
p^2g=pp_{g}+pg_{e}=p^3+pg_{e}\cr
pg_{e}+pg_{p}=pg^2=p_{g}g+g_{e}g=p_{g}g+g_{s}=p^2g+g_{s}
}
$$
et en ajoutant les extr\'emit\'es gauches et droites de ces deux 
lignes on a la formule
$$
\hbox{II})\ \ pg_{p}=p^3+g_{s}
$$
et de mani\`ere analogue on trouve encore la formule (voir \cite[page 26]{Schubert:1879})~:
$$
\hbox{III} )\ \ pg_{s}=p^2g_{p}=G+p^3g=G+p^2g_{e}
\quad .
$$

\paragraph{Justification par la cohomologie}
Si on regarde $\G$ comme l'espace des $2$-plans de ${\mathbb C}^4$,
$\ps$ est le fibr\'e en projectif associ\'e au fibr\'e tautologique 
$\eta$ de rang $2$ sur $\G$. On peut d\'efinir le fibr\'e en droite 
tautologique $\gamma=(F\stackrel{\pi}{\to}\ps)$ par
$$
F=\ens{(\alpha,\ell,v)\in\G\times\p\times{\mathbb C}^4\mid 
\ell\subset\alpha\; ,\, v\in\ell}
\qq
\pi(\alpha,\ell,v)=(\alpha,\ell)
\quad .
$$
Posons $t=c_{1}(\gamma)$; remarquons que $H^*(\ps)$ est un 
$H^*(\G)$-module via l'homomorphisme induit par la projection 
naturelle $\ps\to\G$. On sait (voir \cite[th\'eor\`eme page 62]{Stong:1968}) que 
l'homomorphisme d'anneaux
$$
H^*(\G)[t]\to H^*(\ps)
\qq
t\mapsto c_{1}(\gamma)
$$
induit un isomorphisme
$$
H^*(\G)[t]/I(t^2-tc_{1}(\eta)+c_{2}(\eta))\stackrel{\simeq}{\longrightarrow}
H^*(\ps)
$$
o\`u $I(t^2-tc_{1}(\eta)+c_{2}(\eta))$ d\'esigne l'id\'eal engendr\'e 
par le polyn\^ome $t^2-tc_{1}(\eta)+c_{2}(\eta)$, qui d'ailleurs 
n'est autre que $c_{2}(\eta/\gamma)$ si on fait la substitution
$t=c_{1}(\gamma)$. Il est nul parce que $\gamma$ est un 
sous-fibr\'e de $\eta$.

Prenons la libert\'e d'exprimer par le m\^eme symbole une condition 
sur les \'el\'ements de base et la classe duale par Poincar\'e du 
cycle que cette condition d\'efinit; par exemple, on pourra \'ecrire, 
d'apr\`es le \S\ 2.2.1,
$g=s_{1}(\eta)$. On \'ecrira aussi $t$ pour $c_{1}(\gamma)$.

Il suit de la d\'efinition m\^eme que 
$pg=(-t)(-c_{1}(\eta))$. D'autre part, $p_{g}=t^2$ et on a vu au 
\S\ 2.2.1 que
$g_{e}=c_{2}(\eta)$. On a la relation
$$
tc_{1}(\eta)=t^2+c_{2}(\eta)
$$
que l'on peut \'ecrire encore~:
$$
pg=p^2+p_{g}
$$
qui est la formule I). Donc cette formule est exactement la relation 
par laquelle il faut quotienter $H^*(\G)[t]$ pour obtenir $H^*(\ps)$.

\section{Les formules de co\"\i{}ncidence}

Si $X\subset{\mathbb P}^1\times{\mathbb P}^1$ est une courbe de 
bidegr\'e $(p,q)$, la restriction de son \'equation \`a la diagonale 
est de degr\'e $p+q$, et donc $X$ rencontre cette diagonale en $p+q$ 
points, compt\'es avec leur multiplicit\'e. On peut 
exprimer ce calcul en disant que $X$ est une famille \`a un 
param\`etre de paires de points sur la droite, et qu'il y a $p+q$ 
paires dont les points 
viennent \`a co\"\i{}ncider; c'est ce qu'on appelle le {\it Principe de 
correspondance}\/ de Chasles \cite[Lemme I, page 1175]{Chasles:1864}. 

On va g\'en\'eraliser cette formule en suivant
\cite[ pages 42 et suivantes]{Schubert:1879}, en respectant la notation 
et son ambigu\"\i{}t\'e. On consid\`ere
les paires de points de l'espace, en imaginant que si deux 
points d'une m\^eme paire viennent \`a co\"\i{}ncider, la droite qui 
les joint est encore bien d\'etermin\'ee. 
On appelle $p$ et $q$ les deux points d'une paire, $g$ la 
droite qui les joint et $\e$ la condition que $p$ 
et $q$ sont infiniment proches, mais d\'eterminent encore la droite 
qui les joint. 

On suppose donn\'e un syt\`eme \`a un param\`etre $X$ 
de telles paires de points. Notons que si on d\'esigne encore par $p$ 
le nombre de paires $(P,Q)\in X$ telles que $P$ soit dans un plan 
donn\'e (condition que l'on note aussi $p$), par $q$ le nombre de 
paires telles que $Q$ soit dans un plan donn\'e (condition que l'on 
note aussi $q$), alors $X$ est de bidegr\'e $(p,q)$.

Maintenant on se donne une droite $\ell$ et on 
consid\`ere les paires de plans qui contiennent $\ell$, telles que le premier 
plan contient $p$ et le deuxi\`eme $q$. Ces paires de plans forment 
une courbe $Y$ dans l'espace des paires de plans par $\ell$ (qui 
s'identifie \`a ${\mathbb P}^1\times{\mathbb P}^1$ ; $Y$ est aussi de 
bidegr\'e $(p,q)$), de sorte 
que, par le principe de correspondance de Chasles, il y a
$$
p+q
$$
plans qui contiennent les deux points $p$ et $q$ d'une m\^eme paire. 
Parmi ces plans, on trouve d'une part les $\e$ plans par 
$\ell$ qui contiennent une paire du syst\`eme $X$, dans laquelle les deux points 
de la paire co\"\i{}ncident, et d'autre part aussi les plans par 
$\ell$ qui contiennent une droite $g$ qui joint les deux points 
distincts d'une paire de $X$ (i.e. la droite qui joint les deux points 
de la paire doit couper une droite fix\'ee, \`a savoir la droite 
$\ell$ -- on note $g$ une telle condition, comme d'habitude). On a donc la formule~:
$$
\e=p+q-g
\quad .
$$
Cette formule se r\'ev\`ele tr\`es utile pour \'etablir des 
formules \'enum\'eratives concernant les positions sp\'eciales de 
droites par rapport \`a une surface, par exemple.
Pour la justifier, consid\'erons l'espace $\p\tilde\times\p$ obtenu en 
\'eclatant la diagonale $\Delta$ dans $\p\times\p$. L'application 
qui \`a $(P,Q)\in\p\times\p\setminus\Delta$ associe la droite par $P$ 
et $Q$ s'\'etend en une application $\f:\pt\to\G$, qui est telle que 
$$
\f^*(\Lambda^2(\eta))=\gamma^*\otimes({\cal O}(1)_{1}\otimes{\cal 
O}(1)_{2})
$$
o\`u $\gamma$ d\'esigne le fibr\'e associ\'e \`a la diagonale 
\'eclat\'ee, ${\cal O}(1)_{i}$ le pull-back par la projection de $\pt$ sur le i-\`eme 
facteur de  $\p\times\p$ du fibr\'e des $1$-formes homog\`enes sur $\p$.
On peut s'en convaincre par exemple en utilisant le plongement de 
Pl\"ucker $\psi$ de la grassmannienne $\G$ dans ${\mathbb P}^5$, qui est 
d\'efini ainsi~: si $g\in\G$, on choisit $P,Q\in g$ distincts; si 
$P=[x_{1},\dots ,x_{4}]$, $Q=[y_{1},\dots ,y_{4}]$, et  $x=(x_{1},\dots 
,x_{4})$, $y=(y_{1},\dots ,y_{4})$, on pose 
$$
\psi(g)=[x\wedge y]\in{\mathbb P}(\Lambda^2({\mathbb C}^4))\simeq{\mathbb 
P}^5
\quad .
$$
On v\'erifie que $\psi$ est bien d\'efini et qu'elle d\'efinit un 
plongement, dont l'image est 
$$
\ens{[P\wedge Q]\in {\mathbb 
P}(\Lambda^2({\mathbb C}^4))\mid P\wedge Q\not=0}
\quad ,
$$ 
et on peut 
l' identifier \`a $\G$; notons que le fibr\'e ${\cal O}(1)_{{\mathbb 
P}(\Lambda^2({\mathbb C}^4))}$ des 1-formes homog\`enes sur 
${\mathbb P}(\Lambda^2({\mathbb C}^4))$ est transport\'e par $\psi$ 
sur $\Lambda^2(\eta^*)$.
Consid\'erons l'application 
$$
\Phi:{\mathbb C}^4\times{\mathbb C}^4\to\Lambda^2({\mathbb C}^4)
\qq
(x,y)\mapsto x\wedge y
\quad .
$$
Sa d\'eriv\'ee par rapport \`a $x$, en un point $(y,y)$, $y\not=0$, 
s'\'ecrit $v\mapsto v\wedge y$, qui a pour noyau la droite par $y$. 
Il suit facilement de l\`a que $\Phi$ induit un morphisme
$$
\f:\pt\to\G
\hbox{\quad avec\quad}
\f^*(\Lambda^2(\eta^*))\simeq
\gamma^*\otimes({\cal O}(1)_{1}\otimes{\cal 
O}(1)_{2})
\quad .
$$
Avec un effort suppl\'ementaire, on peut m\^eme montrer que 
$$
\f^*(\eta)=(\gamma^*\otimes{\cal O}_{1}(-1))\oplus{\cal O}_{2}(-1)
$$

On se souvient maintenant que la clase duale \`a $\Omega_{g}$ est 
$s_{1}(\eta)=-c_{1}(\eta)$~; or 
$$
\f^*(-c_{1}(\eta))=
c_{1}(\gamma^*\otimes{\cal O}(1)_{1}\otimes{\cal O}(1)_{2})
$$
et en posant $t_{i}=c_{1}({\cal O}(1)_{i}$, $i=1,2$, 
$\e=c_{1}(\gamma)$, on obtient~:
$$
\f^*(-c_{1}(\eta))=t_{1}+t_{2}-\e
\quad .
$$
Pour retrouver la formule de co\"\i{}ncidence de Schubert, il faut 
encore remarquer que ce qu'on a not\'e $g$ dans ce contexte 
correspond \`a $\f^*(s_{1}(\eta))$, i.e. la condition que la droite 
par $P$ et $Q$ touche une droite donn\'ee $g$, et que $p$ et $q$ 
correspondent \`a $t_{1}$ et $t_{2}$ 
\paragraph{Calculs de co\"\i{}ncidences des intersections d'une droite
avec une surface donn\'ee}

Soit $F\subset\p$ une surface lisse de degr\'e $n$; suivant 
\cite[page 229]{Schubert:1879}, on
d\'esigne par $p_{1}$, $p_{2}$, \dots $p_{n}$ les points 
d'intersection d'une droite $g$ avec $F$. On note $\e_{2}$ la 
condition que $2$ de ces points co\"\i{}ncident. Alors, il suit de 
la formule de co\"\i{}ncidence que
$$
\e_{2}=p_{1}+p_{2}-g
$$
et en multipliant par $g_{s}$~:
$$
\e_{2}g_{s}=p_{1}g_{s}+p_{2}g_{s}-G
$$
en utilisant la formule III)~:
$$
\e_{2}g_{s}=G+p_{1}^3g+G+p_{2}^3g-G=G
$$
car $p^3=0$ (l'intersection g\'en\'erique de $3$ plans et une surface 
est vide). Il reste \`a interpr\'eter le symbole $G$ dans ce 
contexte~: ce sont les paires de points distincts sur l'intersection d'une droite 
g\'en\'erique et la surface, ce qui donne $n(n-1)$ paires. On retrouve 
la formule de la classe d'une courbe plane lisse de degr\'e $n$; en 
effet, $\e_{2}g_{s}$ repr\'esente les droites tangentes \`a la 
surface appartenant \`a un pinceau donn\'e, ce qui revient au m\^eme 
que de consid\'erer les droites du pinceau tangentes \`a la courbe intersection du plan du 
pinceau avec la surface, c'est-\`a-dire les droites par le sommet du 
pinceau tangentes \`a cette courbe.
\paragraph{Justification avec la cohomologie}

D\'esignons 
par $\ft$ l'\'eclatement de $F\times F$ le long de la 
diagonale. On a que $\ft\subset\pt$, et on aimerait conna\^\i{}tre sa 
classe duale; le r\'esultat suivant va nous y aider. Si $Z$ est une 
vari\'et\'e lisse, on notera par $TZ$ son fibr\'e tangent, et par 
$TZ_{x}$ sa fibre au-dessus de $x\in Z$.
\begin{prop}
    Soit $X$ une vari\'et\'e lisse, $A$, $Y\subset X$ des 
    sous-vari\'et\'es lisses, telles que $A\cap Y$ est lisse et que  
    pour tout $x\in A\cap Y$, 
    $TA_{x}\cap TY_{x}=T(A\cap Y)_{x}$, de sorte 
    que l'on a une suite exacte de fibr\'es~:
    $$
    0\to T(A\cap Y)\to TA|_{A\cap Y}\oplus TY|_{A\cap Y}\to
    TX|_{A\cap Y}\to E\to 0
    $$
    o\`u $E$ est  d\'efini par cette suite exacte ;
    on l'appelle {\it fibr\'e exc\`es} et on note $k$ 
    son rang. Alors, en d\'esignant par~:
    
    \vspace{\bigskipamount}
    $$
    \begin{tabular}{c c l}
    $\delta_{U,V}$&\quad& la classe duale \`a $U$ dans $V$\\
    $\tilde X$&\quad&l'\'eclat\'e de $X$ le long de $Y$\\
    $\tilde A$&\quad& le transform\'e strict de $A$\\
    $\e$&\quad&la classe duale au diviseur exceptionnel dans $\tilde X$\\
    $p:\tilde X\to X$&\quad& la projection de l'\'eclatement\\
    $j:\tilde Y\subset \tilde X$&\quad&l'inclusion naturelle
    \end{tabular}
    $$
  \vspace{\bigskipamount}\noindent  
  on a:
    $$
    \delta_{\tilde A,\tilde X}
    =p^*(\delta_{A,X})-\cdot j_{!}\Big((p|\tilde Y)^*(\delta_{A\cap Y,Y})\cdot 
  \underbrace{ 
  \sum_{i=0}^{k-1}(-1)^i\e^ic_{k-i-1}(E)}_{=c_{k-1}(E/\gamma)}\Big)
    $$
\end{prop}
C'est un cas particulier de 
\cite[Theorem 6.7]{Fulton:1984}.

Comme application, on consid\`ere les sous-vari\'et\'es $F\times F$ 
et $\Delta$ de $\p\times\p$. Dans ce cas, le fibr\'e exc\`es 
s'identifie au fibr\'e normal \`a $F$ dans $\p$, qui est ${\cal 
O}(n)_{\Delta}$, et $\delta_{F,\p}=nt$, donc $\delta_{F\times 
F,\p\times \p}=nt_{1}\cdot nt_{2}=n^2t_{1}t_{2}$. Il en suit que
\begin{equation*}
\delta_{\ft,\pt}=n^2t_{1}t_{2}-nt\e\tag{$\heartsuit$}
\end{equation*}
o\`u $t$ d\'esigne indiff\'eremment $t_{1}$ ou $t_{2}$, puisque 
$\e t_{1}=\e t_{2}$.
On peut aussi se convaincre de cette \'egalit\'e de la fa\c con 
suivante; on d\'esigne par
$p:\p\tilde\times\p\to\p\times\p$ la projection 
d'\'eclatement et $\tilde\Delta_{F}=\tilde\Delta\cap^{-1}(F\times 
F\cap\Delta$. Alors~:
$$
p^{-1}(F\times F)=F\tilde\times F\cup\tilde\Delta_{F}
\quad .
$$
En passant aux classe duales, on voit que
$$
\delta_{F\times F,\p\times\p
}=\delta_{\ft,\pt}+\delta_{\tilde\Delta_{F},\pt}
$$
et 
$$
\delta_{\tilde\Delta_{F},\pt}=
\delta_{\tilde\Delta_{F},\tilde\Delta}\cdot\delta_{\tilde\Delta,\pt}
=(nt)\cdot \e
$$
d'o\`u la formule $\heartsuit$.

En particulier, en prenant $n=1$, i.e. $F$ est un plan, on obtient que
$$
\delta_{\f^{-1}(\O_{e})}=t_{1}t_{2}-t\e
\quad .
$$
Or $g_{s}=gg_{e}$, donc 
$\f^*(g_{s})=\f^*(g)\f^*(g_{e})=(t_{1}+t_{2}-\e)(t_{1}t_{2}-t\e)$. 
Pour calculer $\e_{2}g_{s}$ on doit encore multiplier $\f^*(g_{s})\e$ par 
$\delta_{\ft,\pt}$ et \'evaluer cette classe sur $\pt$, ce qui revient 
\`a \'evaluer $\delta_{\ft,\pt}\f^*(g_{s})$ sur $\tilde \Delta$; or
$$
\langle \delta_{\ft,\pt}\cdot\f^*(g_{s}),\tilde\Delta\rangle=
\langle (n^2t_{1}t_{2}-nt\e)(t_{1}t_{2}t\e)(t_{1}t_{2}-\e),\tilde\Delta\rangle=
\langle (n^2t^2-nt\e)(t^2-t\e)(2t-2),\tilde\Delta\rangle
$$
et comme $t^4=0$,
$(n^2t^2-nt\e)(t^2-t\e)(2t-2)=t^2(-n\e^3+\e^2(n^2t+3nt))$. Au lieu 
d'\'evaluer sur $\tilde\Delta$, on peut appliquer $\pi_{!}$ et 
\'evaluer sur $\p$, o\`u $\pi:\tilde\Delta\to\p$ est la projection 
naturelle, qui n'est autre que la projection du fibr\'e en projective 
associ\'e \`a $T\p$; on a ~:
$$
\pi_{!}(\e^2)=1 
\qq
\pi_{!}(\e^3)=c_{1}(T\p)=4t
$$
(cela suit par exemple de la d\'efinition m\^eme des classes de Segre
dans \cite[\S\ 3.1]{Fulton:1984})
et donc
$$
\langle t^2(-n\e+\e^2(n^2t+3nt)),\tilde\Delta\rangle
=\langle t^3(-4n+n^2+3n),\p\rangle=n(n-1)
\quad .
$$

\paragraph{Droites bitangentes} Voici un dernier exemple, que nous 
traiterons \`a la Schubert uniquement (\cite[page 229]{Schubert:1879}).
On note par $\e_{22}$ la condition qu'une droite est tangente en deux 
points de la surface $F$ ; c'est donc la condition que, 
parmi les points $p_1,\dots ,p_n$, 
intersection de la droite avec $F$, deux paires co\"\i{}ncident, 
disons $p_1,p_2$ et $p_3,p_4$. Il suit du principe de co\"\i{}ncidence que~:
$$
2\cdot\e_{22}=(p_1+p_2-g)(p_3+p_4-g)
$$
(le coefficient $2$ venant du fait que l'on peut \'echanger le r\^ole 
de $(p_{1},p_{2})$ et $(p_{3},p_{4})$),
soit
$$
2\cdot\e_{22}=p_1p_3+p_1p_4+p_2p_3+p_2p_4-gp_1-gp_2-gp_3-gp_4+
\underbrace{g^2}_{=g_e+g_p}
\quad .
$$
Les symboles $p_ip_j$, $i\not=j$, ont tous la m\^eme signification, 
de m\^eme que les symboles $gp_i$; on peut donc \'ecrire~:
$$
2\cdot\e_{22}=4p_1p_3-4gp_1+g_e+g_p
\quad .
$$
On multiplie maintenant par $g_e$ : $\e_{22}g_e$ d\'esigne les 
droites bitangentes \`a la surface, situ\'ees dans un plan donn\'e; 
il s'agit donc des droites bitangentes \`a la courbe plane 
intersection de la surface avec le plan. On a~:
$$
2\cdot \e_{22}g_e=4p_1p_3g_e-4p_1gg_e+ g_pg_e=
4p_1p_3g_e-4p_1g_s+G\stackrel{\hbox{(par III))}}{=}
4p_1p_3g_e-4p_1^3g-3G
\quad .
$$
D\'eterminons $p_1p_3g_e$. En fait, nous sommes en train de travailler dans 
$$
\big(F\tilde\times F\big)\times\big(F\tilde\times F\big)\times\G
$$
dont un \'el\'ement g\'en\'erique peut \^etre repr\'esent\'e par 
$((P_1,P_2),(Q_1,Q_2),g)$, avec $P_i,Q_i\in F\cap g$. La condition $p_i$ 
demande que $P_i$ soit dans un plan $e_i$, $i=1,3$, et $g_e$ demande 
que $g$ soit dans un plan $e$. Or $e\cap e_i$ coupe $F$ en $n$ points, 
$i=1,3$; donc la droite des configurations satisfaisant $p_1p_3g$ est 
d\'etermin\'ee par l'une des $n^2$ paires de points, $P_1$ sur 
$e\cap e_1\cap F$ et $P_3$ sur $e\cap e_3\cap F$; pour un tel choix 
de $P_1$ et $P_3$, on peut encore choisir $P_2$ et $P_4$ parmi les 
$n-2$ points restants sur la droite. Il y a donc
$$
n^2(n-2)(n-3)
$$
configurations possibles.
Pour d\'eterminer $G$, il faut remarquer que pour une droite $g$ fix\'ee, 
il y a en tout $n(n-1)(n-2)(n-3)$ paires de points distincts dans $g\cap F$.
Enfin, $p_1^3=0$. On obtient alors~
$$
2\e_{22}g_e=4n^2(n-2)(n-3)-3n(n-1)(n-2)(n-3)
$$
et donc
$$
\e_{22}g_e=\frac12n(n-2)(n-3)(n+3)
\quad .
$$
\bibliographystyle{plain}
\bibliography{schubert}
\end{document}